\documentclass{llncs}

\usepackage[dvips]{graphicx}
\usepackage[latin1]{inputenc}
\usepackage{amssymb,amsmath,array}

\usepackage{epsfig}

\title{Consumer Profile Identification and Allocation}
\author{Patrick Letr\'emy\inst{1}, Marie Cottrell\inst{1}, \\
Eric Esposito\inst{2}, Val\'erie Laffite\inst{2}\and Sally
Showk\inst{2}}

\institute{SAMOS-MATISSE, Universit\'e Paris1-Panth\'eon-Sorbonne,
CES UMR CNRS,\\  90, rue de Tolbiac, F-75013
Paris, France\\
\email{marie.cottrell, patrick.letremy@univ-paris1.fr} \and Research and Development Division\\
Gaz de France\\
\email{eric.esposito, valerie.lafitte,
sally.showk@gazdefrance.com}}

\begin{document}
\maketitle
\begin{abstract}
We propose an easy-to-use methodology to allocate one of the
groups which have been previously built from a complete learning
data base, to new individuals. The learning data base contains
continuous and categorical variables for each individual. The
groups (clusters) are built by using only the continuous variables
and described with the help of the categorical ones. For the new
individuals, only the categorical variables are available, and it
is necessary to define a model which computes the probabilities to
belong to each of the clusters, by using only the categorical
variables. Then this model provides a decision rule to assign the
new individuals and gives an efficient tool to decision-makers.

This tool is shown to be very efficient for customers allocation
in consumer clusters for marketing purposes, for example.
\end{abstract}

\section{Introduction}
The methodology that we propose in this paper is very general and
can be used in many different frames, even if the main
applications belong to the marketing domain. A first presentation
of the main ideas can be found in \cite{Letre}.

Let us define some general notations: Let $\textbf{X}$ be a
database, represented by a $N \times (p+l)$-matrix, where $N$ is
the number of individuals, $p$ the number of continuous variables
(possibly with missing data)  and $l$ the number of categorical
variables (no missing data are allowed). The first $p$ variables
are denoted $X_{1}, X_{2}, \ldots, X_{p}$, and the other $l$
variables are denoted $Y_{1}, Y_{2}, \ldots, Y_{l}$.

In addition, we have a $n \times l$ matrix which corresponds to
new individuals. For these new individuals, only the categorical
variables $Y_{1}, Y_{2}, \ldots, Y_{l}$ are well-informed and the
continuous variables are not available. For example, the $N$
individuals in the database are consumers already registered as
customers by a firm and who can be described by their expenses,
while the new $n$ individuals have only filled up a form and given
some categorical indications (age, housing status, education
level, etc.)
\begin{figure}
\begin{center}
\includegraphics[scale=0.3]{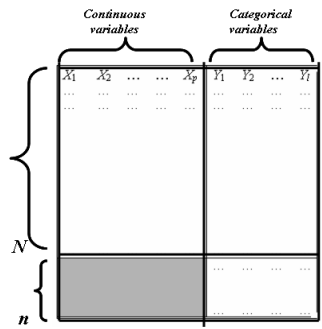}
\caption{The data}
\end{center}
\end{figure}
The first step of the study is to define homogeneous groups from
the point of view of the continuous variables. The interest of
such clustering is double : each cluster corresponds to a typical
profile which is a summary of the whole class and the whole group
can be treated in the same way by any decision-maker. To follow
our example, the direction of sales can use particular targeting
techniques to improve the efficiency of the advertising policy
towards each cluster.

The second step consists in allocating a cluster to new
individuals (for marketing purposes, for example). For this goal,
it is necessary to define a model which computes the probabilities
of belonging to each of the clusters, by using only the
categorical variables. The parameters of this model will be
estimated from the database $\textbf{X}$. For this step, the new
data can be incomplete and missing values are acceptable.

The paper is organized as follows : Section 2 clarifies the
relations between the two types of variables in $\textbf{X}$ and
gives indications about the selection of the relevant variables.
Section 3 briefly deals with the construction of the clusters and
the interpretation of each profile. In Section 4, we present the
multinomial logit model to compute the membership probabilities.
Then this model can be applied to new individuals in order to
assign them to the most probable cluster. Finally section 5 is
devoted to a real-world example and applies the proposed
methodology to a survey data which contains the consumption
structure of Canadian consumers, together with some personal
categorical variables. Section 6 is a conclusion.

\section{Variables selection}
As the final goal is to assign an individual described by the
categorical variables $Y_{1}, Y_{2}, \ldots, Y_{l}$ to a cluster
built from the continuous variables $X_{1}, X_{2}, \ldots, X_{p}$,
it is obvious that the goal cannot be achieved if these two groups
of variables are independent!

We assume that all the categorical variables are of interest for
the applications since they are the only real characteristics
which are available for the new individuals and have to be taken
into account by the decision-makers. So it is necessary to select
the relevant continuous variables which are strongly related to
the categorical ones.

Let us consider the multidimensional $l$-ways additive ANOVA model
( \cite{Rao}) where the explained variables are the $X_{i}$, when
the explanatory variables are the Indicator Functions of each
modality for all the $Y_{j}$. For each component $i, i=1, \ldots,
p$, a global Fisher Statistics and a Squared Correlation
Coefficient are computed. The variables $X_{i}$ that give the
least significant values are not considered in what follows.

\section{The clustering}
First we only take into account the continuous variables $X_{1},
X_{2}, \ldots, X_{p}$ to cluster the $N$ individuals into $K$
clusters. A that step, any unsupervised classification algorithm
can be used. We propose to use a Kohonen algorithm due to several
of its properties (see \cite{Koho1}, \cite{Koho2}, \cite{Kask},
\cite{Cot1}, \cite{Cot2}, \cite{Oja}):

\begin{itemize}
\item The Kohonen maps are known to produce well-balanced and
homogeneous classes,with small quantization error, see
\cite{Verle};
\item The visualization of the clusters is easy to interpret,
thanks to the self-organization property, since there exists a
neighborhood structure between classes;
\item The Kohonen algorithm is robust with respect to missing
values, since it can be adapted to be used with incomplete data,
see \cite{Ibbou}, \cite{Cot3};
\item It is possible to build a Kohonen map having a large number of
classes and to reduce this number by using another clustering of
the code-vectors and thus get a few clusters which will easily be
interpreted and analyzed. These "macro" clusters are composed of
contiguous Kohonen classes, which corroborates the
self-organization property see \cite{Cot2}. To build this second
classification, several methods are available : one can choose an
ascending hierarchical classification or a one-dimensional Kohonen
algorithm. The advantage of this latter choice is that the "macro"
clusters are naturally ordered, and this fact facilitates the
interpretation and the description.
\end{itemize}
After the clusters are built and summarized by their code-vector
or profile, one can describe them from two points of view :
\begin{itemize}
\item The classical statistics (mean, variance, quartile, median)
are computed to characterize and distinguish the clusters.
\item The repartitions of the modalities for each categorical
variable are computed as well as the test values (that is the
ratio between the modality percentage inside the cluster and the
modality percentage in the global population).
\end{itemize}
\section{The model of allocation}
Once the classification stage is achieved, it is necessary to
classify new individuals, who do not belong to the learning set,
in one of the $K$ clusters. A rough method could be to look after
the cluster which contains the number of similar individuals. It
would be a deterministic allocation method. We prefer a stochastic
allocation.

Then, one has to estimate the probability for a new individual to
belong to a cluster only from the categorical variables. The
chosen model is a non-ordered polychotomous logit model, since the
variable to explain (membership probability to a cluster) has more
than two values (there are more than two classes!)

Using the non-ordered polychotomous logit model as discriminating
tool has been proposed by Schmidt and Strauss, \cite{ccc}. It is
an extension of the binary logit model, which is often used in the
studies of appetence or attrition. The explanatory variables are
categorical and the variable to explain can take more than two
modalities \textbf{which are not naturally ordered}. The model
uses the same theoretical frame, since it is estimated by using
the maximum likelihood principle, \cite{bbb}, \cite{ddd}. The
CATMOD procedure of SAS software is designed to estimate this kind
of model, \cite{eee}.

One has to choose a class as a reference class, let us suppose
that it is the class $K$. Let us define $p_{k}=P(k/y)$ as the
probability that an individual belongs to class $k$, given the
fact that it is described by $y$.  Then the non-ordered
polychotomous logit model is written as
$$\frac{p_{k}}{p_{K}}=\exp (y \cdot \beta_k)$$
for $k = 1, 2, \ldots, K-1$, where the $\beta_k \in R^l$ are the
model parameters.

For each possible $y$, the CATMOD procedure provides the estimates
of parameters $\beta_k \in R^p, k=1, \ldots, K-1$ and we compute
the probabilities $p_1, p_2, \ldots, p_{K}$ as a function of $y$
and using equation $p_1 + p_2 + \ldots + p_K=1.$

Then for each new individual $m$ described by $y_m=(y_{jm}), j=1,
\ldots, l$ and for each class $k$, one computes the probability
that this individual belongs to class $k$. The new individual is
assigned to the class for which the probability is maximum.

\section{Application to Canadian consumers data}
We apply the proposed method to a real-world problem: the domestic
consumption of Canadian families. The data have been provided by
Prof. Simon Langlois from the Université of Laval. For 8809
Canadian consumers in 1992, a survey provides the consumption
structures, expressed as percentages of the total expenditure of
the household. Besides, each individual of the survey is also
described by categorical variables (such as Age, Education level,
Wealth, and so on, see below the full list). In a previous
publication, (\cite{Debo}) we have studied the consumptions
profiles, but the allocation problem was not dealt with.

The first step is therefore to identify typical profiles and to
define clusters in the population, on the sole basis of the
consumption structures. Once these profiles and these clusters are
defined, the problem consists in allocating a cluster to new
consumers (for marketing purposes) by using the categorical
variables.

The consumption structure is known through a 19 functions
nomenclature:
\textbf{Consumption nomenclature}. \\
Alcohol; Food at home; Food away; House costs; Communication;
Financial costs; Gifts; Education; Clothes; Housing expenses;
Leisure; Furniture; Health; Security; Personal Care; Tobacco;
Individual Transportation; Collective Transportation; Vehicles.

See in Fig. 2 the mean consumption structure for the 1992 survey.
\begin{figure}
\begin{center}
\includegraphics[width=10cm,height=4cm]{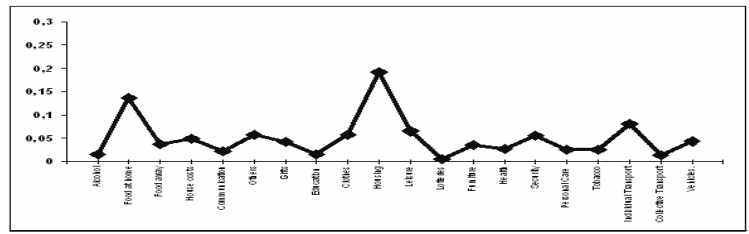}
\caption{Mean Consumption Profile in 1992.}
\end{center}
\end{figure}
For each household, the survey provides also 10 categorical
variables, which concern the head of the family:

\textbf{Categorical variables.} \\
For each item, the number between parenthesis indicates the number
of modalities: Age (4); Language (3); Income (4); Job status (3);
 Professional category (5); Education level (5); Type of town (3);
 Region (5); Residency status (5); Wealth index (5).

We follow the successive steps as described above. First, we write
down the multivariate $l$-ways additive ANOVA model. For 5
consumptions variables (Alcohol; Financial costs; Furniture;
Personal Care; Vehicles), the Squared Correlation Coefficients are
less than $8\%$ and we decide to skip these variables. So, in the
following, we consider that $p=14$ and the percentages are
computed again with only these 14 consumptions functions.

\section{The classification}
We separate the 8809 households into two sets: a learning set with
8400 households, and a test set with 409 households. The consumers
of this test set will further be assigned to one of the cluster,
by using only the categorical variables, and we will compute the
number of correct classifications as a performance measure.

We build a 5-cluster classification in two different ways:
\begin{itemize}
\item First we consider a Kohonen algorithm using a one-dimensional string and 20
units,the number of which is then reduced to 5 macro-classes, by
using another one-dimensional Kohonen algorithm  with 5 units,
that operates over the 20 code-vectors: Classification $C1$.
\item Secondly, more simply, we consider a Kohonen algorithm using a one-dimensional
string with 5 units: Classification $C2$.
\end{itemize}
Fig. 3 represents the 20 Kohonen classes and their code vectors as
well as the 5 macro-classes (marked by different grey). We note
that, due to the topological conservation property of the Kohonen
algorithm, the macro-classes group only neighboring Kohonen
classes. We see that the classes are homogeneous and
well-balanced.
\begin{figure}
\begin{center}
\includegraphics[scale=0.3]{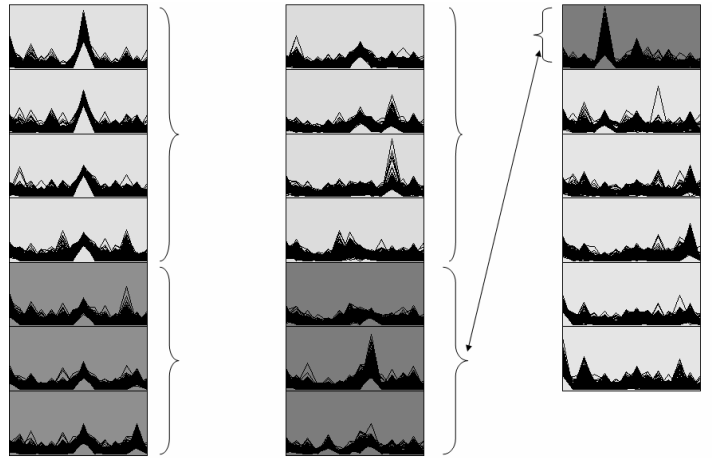}
\caption{The 20 classes from top to bottom and from left to
right and the 5 clusters.}
\end{center}
\end{figure}
Fig 4 shows the code vectors of the 5 macro-classes for the first
classification $C1$.
\begin{figure}
\includegraphics[width=1.0\textwidth]{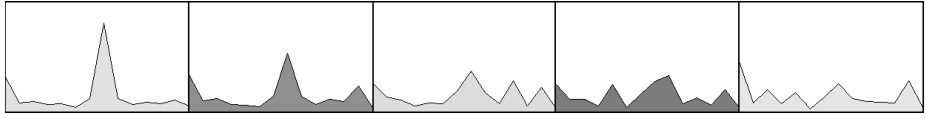}
\caption{The 5 code-vectors for classification $C1$.}
\end{figure}
Fig. 5 shows the 5 code-vectors obtained by
using a one-dimensional Kohonen algorithm with 5 units, classification $C2$.
\begin{figure}
\includegraphics[width=1.0\textwidth]{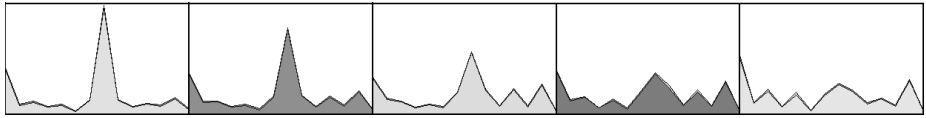}
\caption{The 5 code-vectors for classification $C2$.}
\end{figure}
We see that the code-vectors of $C1$ and $C2$ are similar, and
that the clusters are more or less ordered according to the
housing expenses.

We can represent the distribution of the 4 levels of income across
the 5 classes, in the classification $C2$, see Fig 6.
\begin{figure}
\includegraphics[width=1.0\textwidth]{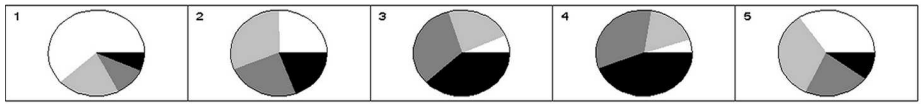}
\caption{Distribution of the 4 income levels over the 5 classes of
$C2$, 4 quartiles from white (low) to black (high).}
\end{figure}
In the same way, we study the distribution of the residency
status, which has 5 modalities: owner without mortgage, owner with
mortgage, tenant, future owner tenant, owner becoming tenant. See
Fig 7.
\begin{figure}
\includegraphics[width=1.0\textwidth]{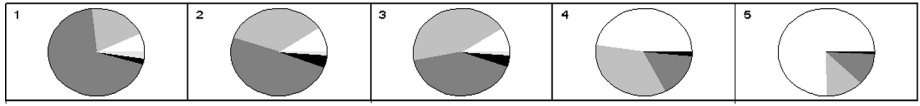}
\caption{Distribution of the 5 residency status modalities over
the 5 classes of $C2$, tenant (medium grey) have majority in class
1, while owner without and with mortgage (white and light grey)
are predominant in classes 3, 4 and 5 .}
\end{figure}
Briefly, we can describe the 5 clusters of $C2$ from the left to
the right. Cluster 1 gathers quasi-poor tenants with low income
and the head of family is unemployed. Cluster 2 groups future
owner tenant, with part time job. In cluster 3, one founds
managers who are owner without mortgage, whore are less than 45
and are quasi-rich or rich. Cluster 4 contains fairly rich workers
who own their house. In cluster 5 there are old  persons who are
owners, but who are poor with a low education level. The housing
expenses are decreasing from cluster 1 to cluster 5, while low
income consumers correspond to clusters 1, 2, 5 and high income
consumers belong to clusters 3 and 4. The tenants are in clusters
1 and 2, the owners are in the others.

In order to evaluate the performance of the allocation algorithm,
we determine:
 1) the "true" class of each individual in the test
set (by comparing its $Y_1, Y_2, \ldots, Y_l$ vector with the
code-vector of each cluster);\\
2) the cluster for which the probability computed in section 4 is
maximum.

We consider the contingency tables that match these two
classifications. The sum of the diagonal entries is the number of
\textbf{exact} allocations. We can also compute the number of
times when the algorithm allocates \textbf{the same cluster or its
neighbors} to the individuals of the test set. This number is the
sum of the entries of the diagonal and of the entries just above
or just below. It is the number of \textbf{correct allocations}.

In the next table, the lines correspond to the
allocation procedure and the columns to the exact classification $C1$.\\

\begin{tabular}{|c|c|c|c|c|c|c|}
  \hline
             & Cluster 1 & Cluster 2 & Cluster  3 & Cluster 4 & Cluster 5 & Total \\
  \hline
  Allocated to 1          & \textbf{55} & \textbf{22} & 29          & 11          &  6    &   123   \\
  Allocated to 2          & \textbf{23} & \textbf{22} & \textbf{14} &  9          &  4    &    72   \\
  Allocated to 3          & 17          & \textbf{11} & \textbf{59} & \textbf{26} &  9    &   122   \\
  Allocated to 4          &  2          &  2          &  \textbf{2} &  \textbf{3} &  \textbf{4} & 13  \\
  Allocated to 5          &  6          &  4          &  7          & \textbf{15} & \textbf{47}  &  79  \\ \hline
  Total & 103  & 61  &  111 &  64  &  70   &   409\\
  \hline
\end{tabular}

From this table, we conclude that the number of exact allocations
is 186, and the number of correct allocations is 186 + 117, that
is 303, which represents $74\%$ of the test set.

We can consider the same table for classification $C2$.\\

\begin{tabular}{|c|c|c|c|c|c|c|}
\hline
               & Cluster 1 & Cluster 2 & Cluster  3 & Cluster 4 & Cluster 5  & Total \\
  \hline
  Allocated to 1          & \textbf{33} & \textbf{12} & 3            & 3            &  5  & 56 \\
  Allocated to 2          & \textbf{23} & \textbf{33} & \textbf{22}  &  17          &  3   &  98 \\
  Allocated to 3          &   8         & \textbf{27} & \textbf{56}  & \textbf{15}  &  3  & 109  \\
  Allocated to 4          &   0         &  3          &  \textbf{11} &  \textbf{42} &  \textbf{21} & 77 \\
  Allocated to 5          &   8         &  3          &  1           & \textbf{10}  & \textbf{47}  & 69 \\
  \hline
  Total & 72  & 78  &  93 &  87 &  79   &   409\\
  \hline
\end{tabular}

The number of exact allocations is 211, and the number of correct
allocations is 186 + 141, that is 352, which represents $86\%$ of
the test set.
\section{Conclusion}
This paper proposes a simple method to study large databases which
contain continuous and categorical variables at the same time, and
to identify the category of an individual who is described by
incomplete data. The results are convincing and it can give a very
useful tool to decision-makers in many fields : insurance
policies, personal tariffing, targeted advertising, credit
scoring, etc.
\section{Acknowledgement}
The authors would like to thank Patrice Gaubert from SAMOS-MATISSE
and Cr\'eteil Universit\'e for making the Canadian consumption
data available to us, and the Gaz de France Company for partially
funding this research via a previous partnership and collaboration
on the allocations problem.

\end{document}